\documentclass{article}
\usepackage{graphicx}
\usepackage{amsmath}
\usepackage{amsfonts}
\usepackage{amssymb}

\begin{document}
Dear prof. Francisco L. Hernandez,

\bigskip

Thank you for your letter.
You are completely right and my paper has a fatal error.
It contains in the last phrases of proof. What a shame on my grey hair!
Moreover, I MUST know that the result is not true because the same ''proof''
may be applied in a separable case too. I sincerely apologize.
Moreover, I knew it in 1972 when was written my first paper (no. [5] in list
of references).
I explain what it was known. May be it will be useful for you in your next work.

\bigskip

It is obvious that if $E$ be a Banach space with a symmetric basis $\left(
e_{n}\right)  $ and $F$ is its subspace with symmetric basis $\left(
f_{n}\right)  $ then $\left(  f_{n}\right)  $ is equivalent to a
block-sequence $\sum_{n_{k}+1}^{n_{k+1}}a_{i}e_{i}$ and there may be two possibilities:

1. $c_{k}=\max\{a_{i}:n_{k}+1\leq i\leq n_{k+1}\}\rightarrow0$

2. $\inf c_{k}\geq c>0$.

Certainly, if $E$ and $F$ are nonseparable (with symmetric uncountable bases)
then the possibility 1 cannot be realized.

So, in nonseparable case the embedding $F\hookrightarrow E$ implies the
continuous inclusion $F\subset^{C}E$, (i.e. if $f=\sum\xi_{k}f_{k}\in F$;
$e=\sum\xi_{k}e_{k}\in E$ then $\left\|  f\right\|  _{F}\geq c\left\|
e\right\|  _{E}$).

Since for every space $E(\varkappa)$ with a (may be, uncountable) symmetric
basis $l_{1}\left(  \varkappa\right)  \subset^{C}E\left(  \varkappa\right)
\subset^{C}c_{0}\left(  \varkappa\right)  $ ($\varkappa\geq\omega$ is a
cardinal), immediately follows the S. Trojanski's result on nonseparable
$l_{1}\left(  \varkappa\right)  $. Indeed, if $E\left(  \varkappa\right)  $
contains a subspace isomorphic to $l_{1}\left(  \varkappa\right)  $ then
$l_{1}\left(  \varkappa\right)  \subset^{C}E\left(  \varkappa\right)
\subset^{C}l_{1}\left(  \varkappa\right)  $ i.e $l_{1}\left(  \varkappa
\right)  \ $is isomorphic to $E\left(  \varkappa\right)  $.

However concerning other $p$'s NOTHING is follows since these spaces has no
any special properties in the partially ordered set (by the order $\subset
^{C}$) of all spaces with symmetric bases of given cardinality $\varkappa$.

This simple observation makes clear reasons by which the J. Lindenstrauss
result (on embeddings into spaces with symmetric bases) cannot be generalized.
Indeed, it is clear that the interpolation construction transforms a basis
$\left(  y_{n}\right)  $ of an embedable space $Y$ into a block-basis that
satisfies the possibility 1, mentioned above (what is impossible, as was
noted, in a nonseparable case).

\bigskip

I'l try to remove my misfortune paper from arXive.

\bigskip

Thank you once more.

I did not see almost any paper since 1990 till today. This partially explain my fault.

\bigskip

I would be very pleasant for your next remarks, comments, suggestions and questions.

\bigskip

Sincerely your

\bigskip

Eugene Tokarev
\end{document}